\documentclass[10pt]{article}
\input{epsf}
\input{amssym}
\def\Cchi{{\raisebox{.2ex}{ \large $\chi$}}}
\newtheorem{prethm}{{\bf Theorem}}

\newenvironment{thm}{\begin{prethm}{\hspace{-0.5
em}{\bf.}}}{\end{prethm}}
\newtheorem{precor}{{\bf Corollary}}

\newtheorem{pretheorem}{{\bf Theorem}}

\newenvironment{theorem}{\begin{pretheorem}{\hspace{-0.85
em}{\bf.}}}{\end{pretheorem}}
\newtheorem{preprop}{{\bf Proposition}}

\newtheorem{preque}{{\bf Question}}

\newtheorem{preques}{{\bf Question}}

\newenvironment{ques}{\begin{preques}{\hspace{-0.85
em}{\bf.}}}{\end{preques}}
\newtheorem{prelemma}{{\bf Lemma}}

\newtheorem{prelemm}{{\bf Lemma}}

\newtheorem{preex}{{\bf Example}}

\newenvironment{ex}{\begin{preex}{\hspace{-0.5
em}{\bf.}}}{\end{preex}}
\newtheorem{prepro}{{\bf Proposition}}

\newtheorem{prelem}{{\bf Theorem}}

\newenvironment{lem}{\begin{prelem}{\hspace{-0.5
em}{\bf}}}{\end{prelem}}
\newtheorem{preproof}{{\bf Proof.}}

\newenvironment{proof}[1]{\begin{preproof}{\rm
               #1}\hfill{$\rule{2mm}{3mm}$}}{\end{preproof}}

\newtheorem{preconj}{{\bf Conjecture}}

\newtheorem{predefinition}{Definition}

\newenvironment{definition}[1]{\begin{predefinition}{\rm
               #1}{}}{\end{predefinition}}
\newtheorem{preremark}{Remark}

\newenvironment{remark}[1]{\begin{preremark}{\hspace{-0.5
em}{\bf.}}{\rm #1}}{\end{preremark}}
\def\newpic#1{}
\date{}
\begin{document}

\title{
{\bf Constructing regular graphs with smallest defining number}
 }
\author{
{\sc Behnaz Omoomi}$^a$ and {\sc Nasrin Soltankhah}$^{b}$
\footnote{This research was in part supported by a grant from the
Institute for Studies in Theoretical Physics and
Mathematics (IPM).}\\
[5mm]
$^a$ Department of Mathematical Sciences \\
Isfahan University  of Technology \\
Isfahan, 84156-83111 \\ [5mm]
$^b$ Department of Mathematics,  Alzahra University\\
Vanak Square 19834,  Tehran, Iran \\
[5mm]
}
\maketitle
\begin{abstract}
In a given graph $G$,  a set $S$ of vertices with an assignment
of colors is a {\sf defining set of the vertex coloring of $G$},
if there exists a unique extension of the colors of $S$ to a
$\Cchi(G)$-coloring of the vertices of $G$. A defining set with
minimum cardinality is called a {\sf smallest defining set\/} (of
vertex coloring) and its cardinality, the {\sf defining number},
is denoted by $d(G,  \Cchi)$.  Let $ d( n, r, \Cchi = k)$
 be the smallest defining number of all $r$-regular $k$-chromatic
graphs with $n$ vertices. Mahmoodian et. al~\cite{rkgraph} proved
that, for a given $k$ and  for all $n \geq 3k$, if  $r \geq
2(k-1)$ then
 $ d( n,  r,  \Cchi = k)=k-1$.
 In this paper we show that
 for a given $k$ and for all $n < 3k$ and $r\geq 2(k-1)$,
$d(n, r, \Cchi=k)=k-1$.
\end{abstract}
\begin{center}
\end{center}
%
\section{Introduction} 
A {\sf k-coloring} of a graph $G$ is an assignment of $k$
different colors to the vertices of $G$ such that no two adjacent
vertices receive the same color. The ({\sf vertex}) {\sf
chromatic number\/}, $\Cchi(G)$, of a graph $G$ is the minimum
number $k$  for which there exists a $k$-coloring for $G$. A
graph $G$ with $\Cchi(G) = k$ is called a {\sf $k$--chromatic}
graph. In a given graph $G$,  a set of vertices $S$ with an
assignment of colors is called a {\sf defining set of vertex
coloring},  if there exists a unique extension of the colors of
$S$ to a $ \Cchi(G) $-coloring of the vertices of $G$. A defining
set with minimum cardinality is called a {\sf smallest defining
set\/} (of a vertex coloring) and its cardinality is the {\sf
defining number},  denoted by $d(G,  \Cchi)$.

There are some results on defining numbers in~\cite{mnz} (see
also {{\cite{k96}}, and {\cite{kerman}}). Here we study the
following concept. Let $d(n, r, \Cchi=k)$ be the smallest value
of $d(G,\Cchi)$ for all $r$-regular $k$-chromatic graphs with $n$
vertices. Note that for any graph $G$, we have $d(G,\Cchi) \geq
\Cchi(G)-1$, therefore $d(n, r, \Cchi=k)\geq k-1$. By Brooks'
Theorem~\cite{brooks}, if $G$ is a connected  $r$-regular
$k$-chromatic graph  which is not  a complete graph or an odd
cycle, then $k \leq r$. For the case of $r=k$, Mahmoodian and
Mendelsohn~\cite{mahmend} determined the value of $d(n, k,
\Cchi=k)$ for all $k \leq 5$. Mahmoodian and
Soltankhah~\cite{mahsol}  determined this value for $k=6$ and
$k=7$. Also   in~\cite{mahsol}, for each $k$, the value of $d(n,
k, \Cchi=k)$ is determined for  some congruence classes of $n$.
For the case of $k < r$, it is proved in~\cite{mahmend} that, for
each $n$
 and each $r \geq 4$, we have $d(n, r, \Cchi=3)=2$.
The following question is raised in~\cite{mahmend}:
\begin{ques}
Is it true that for every $k$, there exist  $n_0(k)$ and $r_0(k)$,
such that
for all $n\geq n_0(k)$ and $r\geq r_0(k)$ we have $d(n, r, \Cchi=k)=k-1${\rm ?}
\end{ques}
 Mahmoodian et. al.~\cite{rkgraph} proved that
 the answer to this question
is positive and  that, for a given $k$ and  all
$n  \geq 3k$, if
$r \geq 2(k-1)$ then  $d(n, r, \Cchi =k)=k-1$.

We show the above statement for $n <3k$.
 In fact we prove that:
\begin{theorem}
Let $k$ be a positive integer. For each $n <3k$, if $r \geq
2(k-1)$ then $d(n, r, \Cchi =k)=k-1$.
\end{theorem}
%
\section{Preliminaries}   
\indent   In this section, we state some known results and definitions
  which will be used in the sequel.
\begin{definition}{{\rm \cite{mahmend}}.
Let $G$ and $H$ be two graphs,  each with a given proper
$k$-{\mbox coloring}
 say $c_G$ and $c_H$, (respectively) with $k$ colors.
Then the {\sf chromatic join} of $G$ and $H$,  denoted by
$G\stackrel{\chi}{\vee}H$
 is a graph where  $V (G\stackrel{\chi}{\vee}H)$ is  $V(G)\cup V(H)$,
 and $E (G\stackrel{\chi}{\vee}H)$ is  $E(G) \cup E(H)$,  together with
the set $\{ xy \ |\ x\in V(G),\ y\in V(H)\ {\rm such\ that} \ c_G
(x) \neq c_H (y) \}.$
 }\end{definition}
%
\begin{lem} {\rm \cite{mahmend}}. \label{kl}
Let $n$ be a multiple of $k$, say $n=kl$ $(l \geq 2)$; then
$d(kl, 2(k-1), \Cchi=k)=k-1$.
\end{lem}
To prove this  theorem Mahmoodian and Mendelsohn constructed a
$2(k-1)$-regular $k$-chromatic graph with $n=kl$ vertices as
follows. Let $ G_1,  G_2,  \ldots,  G_l  $ be vertex disjoint
graphs such that $G_1$ and $G_l$ are  two copies of $K_k$ and if
$l \geq 3$,  $G_2, \cdots,  G_{l-1}$ are copies of
$\overline{K}_k$. Color each $G_i$ with $k$ colors $1,  2,
\cdots,  k$. Then construct a graph $G$ with $lk$ vertices by
taking the union of $ G_1 \cup G_2 \cup \ldots \cup G_l $, and by
making a chromatic join between $G_i$  and $G_{i+1}$; for $i = 1,
2,  \cdots,  l-1$. This is the desired graph. We denote such a
graph by $G_{l(k)}$ and use  this construction in Section~3.
 \begin{definition} \ {\rm{\cite{mahsol}}.
Let $G$ be a $k$-chromatic graph and let $S$ be a defining set
for~$G$. Then a set $F(S)$ of edges is called {\sf nonessential
edges}, if the chromatic number of  $G - F(S)$, the graph
obtained from $G$ by removing the edges in $F(S)$, is still $k$,
and $S$ is also a defining set for $G-F(S)$.
 }\end{definition}
\begin{remark} {\rm \ A necessary condition for the existence of an $r$-regular $k$-chromatic
graph is $\frac{r}{k-1} \leq \frac{n}{k}$. For,  if $G$ is an
$r$-regular $k$-chromatic graph with $n$ vertices, then each
chromatic class in $G$ has at most $n-r$ vertices. Therefore $n
\leq k(n-r)$. This implies  $\frac{r}{k-1} \leq \frac{n}{k}$.
Thus, for $r \geq 2(k-1)$ there are not any graph of order $n <
2k$. Hence when $r \geq 2(k-1)$, it is sufficient to investigate
$d(n, r, \Cchi=k)$ only for $n \geq 2k$. Also it is obvious that
$n$ and $r$ cannot be both odd. }\end{remark}
For the  definitions and  notations not defined here we refer
the reader to texts,
such as \cite{west}.
%
\section{Main results}      
\indent In this section in the following four theorems we prove
our main result, which was mentioned at the end of Section~1.
\begin{thm}
For each $k \geq 3$ and each $r \geq 2(k-1)$, we have \\
$d(3k-1,r, \Cchi=k)=k-1$.
\end{thm}
\begin{proof}{
Let  $n=3k-1$ and $r=2(k-1)+t$. By Remark~1 it is obvious that $t
\leq k-2$. First for $t=0$, we construct a $2(k-1)$-regular
$k$-chromatic graph $H$ with $n$ vertices and $d(H, \Cchi)=k-1$
as follows. By Theorem A  we have $d(3k, 2(k-1), \Cchi=k)=k-1$.
In   graph $G_{3(l)}$ which was constructed to prove Theorem A,
let $V(G_1)=\{ u_1, u_2,...,u_k \}$, $V(G_2)=\{  v_1, v_2,...,v_k
\}$, and $V(G_3)=\{  w_1, w_2,...,w_k \}$. Also assume that
$c(u_i)=c(v_i)=c(w_i)=i$, for $i=1,2,...,k$. Note that the set
of  vertices adjacent to $v_k$ is $N_{G_{3(l)}} (v_k)=\{
u_1,...,u_{k-1} \} \cup \{  w_1,...,w_{k-1} \}$. We delete the
vertex $v_k$ and join its neighbors in the following manner: we
join $u_i$ to $w_{i+1}$ for $i=1,2,...,k-2$ and $u_{k-1}$ to
$w_1$. It can be easily seen that the new graph, say $H$, is
$2(k-1)$-regular $k$-chromatic  with $n=3k-1$ vertices with a
defining set $S=\{  u_1,u_2,...,u_{k-1} \}$.

Now for $1 \leq t \leq k-3$,  to construct an $r$-regular $k$-chromatic
graph, we consider the graph $H$, and we add the edges
$u_i w_{i+j+2}$ $\pmod {k}$, for $i=1,...,k$ and $j=1,...,t$, to $H$.
Also, in the case of $k$ odd, we add the edges of $t$ mutually disjoint
1-factors of $K_{k-1}$, and in the case of $k$ even,
 the edges of $\frac{t}{2}$
mutually disjoint 2-factors of $K_{k-1}$,  on vertex set
$\{v_1,...,v_{k-1}\}$.

Note that if $t=k-2$ then such a graph does not exist.
For, if $G$ is a graph satisfying  such conditions then we know that each chromatic
class in  $G$ has at most $3$
vertices. Since $n=3k-1$, $G$ must have $k-1$ chromatic classes of
size $3$ and one chromatic class of size $2$. And each vertex in
a chromatic class of size $3$ must be adjacent to all other vertices.
This implies that the degree of each vertex in the chromatic class of size
$2$ is $3(k-1)=r+1$, which  contradicts the  $r$-regularity of the graph.
}\end{proof}
\begin{ex}
In Figure~\ref{5.6} we show the graph $H$ when $k=5$ and $r=8$.
 The vertices of the defining set are shown by the filled circles.
\begin{figure}[h]
\unitlength=1mm
 \centering
 \epsfysize=4.8cm
 \epsfbox[0 5 300 130]{G3(5).eps}
\caption{\label{5.6} $d(H, \chi=5)=4$.}
\end{figure}
\end{ex}

\begin{thm} \label{kodd}
For each odd number $k \geq 3$,  and each $2k \leq n \leq 3k-2$,
we have $d(n , 2(k-1), \Cchi=k)=k-1$.
\end {thm}
\begin{proof}{
 By Theorem A  we have $d(2k , 2(k-1), \Cchi=k)=k-1$.
Let $n=2k+s$, $s=1,2,...,k-2$. We construct a $2(k-1)$-regular $k$-chromatic
graph $H_s$ with $n$ vertices and defining number equals to $k-1$.
For this, we consider  graph
 $G_{2(l)}$ and add $s$ new vertices
 to it, delete some suitable edges as follows and join the new vertices
 to the end vertices of deleted edges. In  graph $G_{2(l)}$, for convenience
let $V(G_1)=\{  u_1,...,u_i,...,u_{{k-1} \over {2}},
 u_{1'},...,u_{i'},...,u_{{({{k-1} \over {2}})}'}\\,u_k \}$
 and  $V(G_2)=\{  v_1,...,v_i,...,v_{{k-1} \over {2}},
 v_{1'},...,v_{i'},...,v_{{({{k-1} \over {2}})}'},v_k \}$,
  where  $i'=i+{{k-1} \over {2}}$, $i=1,2,...,{{k-1} \over {2}}$;
 and $c(u_j)=c(v_j)=j$, for $j=1,2,...,k$.

\noindent If $ 1 \leq s \leq {{k-1}\over {2}}$ then denote new
vertices by $x_1,...,x_s$. Let $M_1, M_2,..., M_{{k-1}\over{ 2}}$
be mutually
 disjoint 1-factors of subgraph
 $< u_1,...,u_i,...,u_{{k-1} \over {2}},
 u_{1'},...,u_{i'},...$ $,u_{{({{k-1} \over {2}})}'} >$
in $G_{2(l)}$ such that each edge in $M_i$ has one end in $\{
u_1,u_2,...,u_{{k-1}\over{ 2}} \}$ and the other end in $\{
u_{1'},...,,u_{{({{k-1} \over {2}})}'} \}$. For each $i$ $(1 \leq
i \leq s)$ we join $x_i$ to each of the vertices of $M_i$, and
delete all of the edges of $M_i$. Also with respect to each $u_a
u_b \in M_i$, we delete the edge $v_a v_b$ and join $x_i$ to the
vertices $v_a$ and $v_b$. Now it can be easily seen that
$deg(x_i)=2(k-1)$. Note that the new graph contains a complete
subgraph say,  $< u_1,u_2,...,u_{{k-1}\over {2}},
v_{1'},...,v_{{({{k-1} \over {2}})}'}, x_1>=K_k$ and a defining
set $S=\{u_1,...,\\u_{k-1}\}$. Also the colors of vertices of
$G_{2(l)}$ force  all new vertices to be colored~$k$.

\noindent If ${{k-1}\over{ 2}} < s \leq k-2$ then we denote the
new vertices by $x_1,x_2,...,x_{{k-1}\over {2}}, y_1,\\ y_2,...,
y_{s-{{k-1} \over {2}}}$.
 For $x_i$ $(1 \leq i \leq {{k-1}\over {2}})$ we proceed as before.
For $y_t$ $(1 \leq t \leq s-{{k-1} \over {2}})$, first we
recognize some nonessential edges in $H_{{k-1}\over {2}}$.
If for each $i$, we let $z_i$ be either $u_i$ or $v_i$ and, for each
$j$, we let $w_j$ be either $u_j$ or $v_j$, then the following
edges form a nonessential set in $H_{{k-1}\over {2}}$:

$\begin{array}{cl}
 \hspace{-8mm} F=&\{ v_i v_j   \   |\  1 \leq i
< j \leq {{k-1}\over {2}}   \}     \cup
\{ u_{i'} u_{j'} \  |\  1' \leq i' < j' \leq  {({{k-1}\over{ 2}})}' \} \cup \\
&  \{ x_1 u_{i'} \ or \  x_1 v_i  \ |\ 1 \leq i \leq {{k-1}\over
{2}} \}
 \cup \{ x_i w_j \ |\  2 \leq i \leq {{k-1}\over {2}} , 1 \leq j \leq k-1 \}
 \cup \\ & \{ z_i v_k \ |\  1 \leq i \leq k-1 \}.
 \end{array} $

There are two cases to be considered.\\
[5mm]
 {\bf  Case 1.} $k=4l+1$. \\
[5mm] \indent In this case the induced subgraphs
$A=<u_{1'},u_{2'},...,u_{{({{k-1} \over {2}})}'}>$ and
$B=<v_1,v_2,..., v_{{k-1}\over {2}} >$  are complete graphs
$K_{{k-1} \over {2}}$. So they  are 1-factorable. Let $F_1,
F_2,..., F_{{k-3}\over{ 2}}$ and $F'_1, F'_2,..., F'_{{k-3}\over{
2}}$  be  1-factorizations of $A$ and $B$, respectively, such
that the edge $u_{t'} u_{{({{k-1} \over {2}}) }'} \in F_t$ and
$v_t v_{{k-1} \over {2} } \in F'_t$. Now for each $t$ $(1 \leq t
\leq {s-{{k-1} \over {2}}} \leq {{k-3} \over {2}})$ we delete all
of the edges of $F_t \backslash \{u_{t'} u_{{({{k-1} \over
{2}})}' }  \}$ and $F'_t \backslash \{ v_t v_{{k-1} \over {2}}
\}$. Also we delete the edges $u_t v_{{k-1} \over {2}}$ and
$u_{t'} v_k$. Finally we delete all the edges $x_1 v_t$, $x_2
u_{t+1}$,..., $x_{{k-1} \over 2} u_{t+{{k-3} \over 2 }}$ $\pmod
{{{k-1} \over 2}}$. We join $y_t$ to the ends of all   deleted
edges. It can be easily  seen that $deg(y_t)=2(k-1)$ and the
color of
$y_t$ is forced to be $k-1$.\\
[5mm]
 {\bf Case 2.} $k=4l+3$.\\
[5mm] \indent In this case the induced subgraphs
$A=<u_{1'},u_{2'},...,u_{{({{k-1} \over {2}})}'}, u_k>$ and
$B=<v_1, v_2,...,v_{{k-1}\over {2}}, v_k >$ are complete graphs
$K_{{k+1}\over {2}}$. Thus they are 1-factorable. Let $F_1,
F_2,..., F_{{k-1}\over{ 2}}$ and $F'_1, F'_2,..., F'_{{k-1}\over{
2}}$ be 1-factorizations of $A$ and $B$, respectively, such that
$u_{t'}u_k \in F_t$ and $ v_tv_k \in F'_t$, for $1 \leq t \leq
{{k-1}\over{ 2}}$. Now for each $t$ $(1 \leq t \leq {s-{{k-1}
\over {2}}} \leq {{k-3} \over {2}})$ we delete all of the edges
of  $F_t \backslash \{ u_{t'}u_k \}$  and $F'_t \backslash \{
v_tv_k \}$. Also we delete the edge $v_k u_t$.
 Finally we delete the edges $x_1 v_t$, $x_2 u_{t+1}$,...,
$x_{{k-1} \over 2} u_{t+{{k-3} \over 2 }}$ $\pmod {{{k-1} \over 2}}$.
We join $y_t$ to the ends of all  deleted edges.
It can be easily seen that $deg(y_t)=2(k-1)$ and the  color of $y_t$
 is forced to be $t+{{k-1} \over {2}}$.
}\end{proof}
To illustrate the construction shown  in the  proof of
Theorem~\ref{kodd}, we provide the following example.
\begin{ex}
{\rm Let $k=7$. For $15 \leq n \leq 19$, we construct a
$12$-regular $7$-chromatic graph of order $n$  with a defining
set of size $6$. For $n=14+s$, $1 \le s \le 5$, we add $s$ new
vertices to the $12$-regular $7$-chromatic graph  $G_{2(7)}$ of
order $14$ and delete some nonessential edges as  explained in the
proof of Theorem~\ref{kodd}.

\begin{table}[h]
\begin{center}{
\caption{\label{t5.1} New vertices and deleted edges.}
\vspace{.25cm}
\begin{tabular}{|c|c|c|c|c|c|}
\hline
New vertices  &$x_1$         & $x_2$       & $x_3$      & $y_1$                & $y_2$ \\
\hline
& $u_1u_{1'}$ & $u_1u_{2'}$ & $u_1u_{3'}$ & $u_{2'}u_{3'}$ & $u_{1'}u_{3'}$\\
&$u_2u_{2'}$ & $u_2u_{3'}$ & $u_2u_{1'}$ & $v_2v_3$       & $v_1v_3$       \\
Deleted &$u_3u_{3'}$ & $u_3u_{1'}$ & $u_3u_{2'}$ & $x_1v_1$       & $x_1v_2$       \\
edges &$v_1v_{1'}$ & $v_1v_{2'}$ & $v_1v_{3'}$ & $x_2u_2$       & $x_2u_3$       \\
&$v_2v_{2'}$ & $v_2v_{3'}$ & $v_2v_{1'}$ & $x_3u_3$       & $x_3u_1$       \\
&$v_3v_{3'}$ & $v_3v_{1'}$ & $v_3v_{2'}$ & $v_7u_1$       & $v_7u_2$       \\
\hline
\end{tabular}
}\end{center}
\end{table}
Table~\ref{t5.1} gives all  the deleted edges of $G_{2(7)}$ with
respect to  addition of new vertices. In Figure~\ref{f5.2}, we
show  the deleted edges and the added edges to construct a
$12$-regular $7$-chromatic graph $H_1$ of order $15$ $(s=1)$ with
a defining set of size $6$.  The dotted lines are the deleted
edges and the vertices of the defining set are shown by the
filled circles.
\vspace{4mm}\\
\begin{figure}[ht]
\unitlength=1mm
 \centering
 \epsfysize=6cm
 \epsfbox[60 10 200 120]{G2(7)2.eps}
\caption{\label{f5.2} $d(H_1,\chi=7)=6$.}
\end{figure}
}\end{ex}
\begin{thm}\label{keven}
For each even number $k \geq 4$, and each $2k \leq n \leq 3k-2$,
we have $d(n , 2(k-1), \Cchi=k)=k-1$.
\end {thm}
\begin{proof}{
 By Theorem A  we have $d(2k , 2(k-1), \Cchi=k)=k-1$.
For $n=2k+s$, $s=1,2,...,k-2$, we construct a $2(k-1)$-regular $k$-chromatic
graph $H_s$ with $n$ vertices and defining number equal to $k-1$.\\
 To construct $H_s$, we consider  graph $G_{2(k)}$
and add $s$ new vertices to it, delete some suitable edges and
 join the new vertices to the end vertices of the deleted edges as follows.
  In  graph $G_{2(k)}$ for convenience let
 $V(G_1)= \{ u_1,...,u_i,...,u_{{k} \over {2}}, u_{1'},...,u_{i'},...,
 u_{{({{k} \over {2}})}'} \}$ and
  $V(G_2)= \{ v_1,...,v_i,...,v_{{k} \over {2}}, v_{1'},...\\,v_{i'},...,
v_{{({{k} \over {2}})}'} \}$,  where  $i'=i+{{k} \over {2}}$,
 $i=1,2,...,{{k} \over {2}}$; and $c(u_j)=c(v_j)=j$, for $j=1,2,...,k$.

If $ 1 \leq s \leq {{{k}\over {2}}-1}$ then we denote the new
vertices by $x_1,...,x_s$. Let  $M_1, M_2,...,M_{{k}\over{ 2}}$ be
mutually disjoint 1-factors of the induced subgraph
$G_1=<u_1,...,u_i,...,u_{{k} \over {2}}, u_{1'},...,u_{i'},...,
u_{{({{k} \over {2}})}'} >$, where, for $i=1,2,...,{{k} \over
{2}}; $\\
$$M_i=\{u_1u_{i'},u_2u_{{(i+1)}'},...,u_tu_{{(i+t-1)}'},...,u_{{k}
\over {2}} u_{{(i+{{k} \over {2}}-1)}'} \} \pmod {{ {k} \over {2}
}}. $$
 Also let $M'_1,
M'_2,...,M'_{{k}\over{ 2}}$ be mutually  disjoint 1-factors of
the induced subgraph
 $G_2= <v_1,...,v_i,...,v_{{k} \over {2}},v_{1'},...,v_{i'},...,
 v_{{({{k} \over {2}})}'}>,$ where, for $i=1,2,...,{{k} \over {2}}; $
$$M'_i = \{ v_1 v_{i'},v_2 v_{{(i+1)}'},...,v_t v_{{(i+t-1)}'},..
.,v_{{k} \over {2}} v_{{(i+{{k} \over {2}}-1)}'} \} \pmod {{ {k}
\over {2}} }. $$
 Now for each $i$ $(i=1,2,...,s)$ we delete  all of the edges of
 $M_{i+1} \backslash \{u_{{{k} \over {2}}-i} u_{{({{k} \over {2}})}'} \} $,
 and  all of the edges of
 $M'_i \backslash \{ v_{{{k} \over {2}}-i+1} v_{{({{k} \over {2}})}'} \}$.
 Finally we delete the edge $u_{{{k} \over {2}}-i} v_{{{k} \over {2}}-i+1}$.
We join $x_i$ to the ends of all  deleted edges.
Now it can be easily seen that $deg(x_i)=2(k-1)$.
 Note that the new graph contains a complete subgraph say
 $< u_1,u_2,...,u_{{k}\over {2}},u_{{({{k}\over {2}})}'},
v_{1'},\linebreak[0]...,v_{{({{k} \over {2}}-1)}'}>=K_k$ and a
defining set $S=\{u_1,...,u_{k-1}\}$. Also the colors of vertices
of $G_{2(k)}$ force the colors of all new vertices to be $k$.

 If ${{k}\over{ 2}} \leq  s \leq k-2$ then we denote the new vertices by
$x_1,x_2,...,x_{{{k}\over {2}}-1}, y_1, y_2,\linebreak[0]...,
 y_{s-{{k} \over {2}}+1}$. For $x_i$
$(1 \leq i \leq {{{k}\over {2}}-1})$ we treat as before.
For $y_t$ $(1 \leq t \leq s-{{k} \over {2}}+1)$ first we
recognize some nonessential edges in $H_{{{k}\over {2}}-1}$.
If for each $j$, we let $w_j$ be either $u_j$ or $v_j$,
then the following edges form a
 nonessential set in $H_{{{k}\over {2}}-1}$:

$\begin{array}{cl}
 \hspace{-6mm} F=&\{ v_i v_j \ |\  1 \leq i  < j
\leq {{k}\over {2}}, j \neq {i+1}   \} \cup
  \{ u_{i'} u_{j'} \ |\  1' \leq i'< j' \leq {{{({{k}\over{ 2}})}'}-1} \}
 \cup \\ & \{ x_i w_j\ |\  1 \leq i \leq
{{{k}\over {2}}-1} , 1 \leq j \leq k-1 \} \cup
\{ v_i v_{{({{k} \over {2}})}'} \ |\
1 \leq i \leq {{{({{k} \over {2}})}'}-1} \} \cup \\
 & M_1 \backslash \{u_{{k} \over {2}} u_{{({{k} \over {2}})}'}\}  \cup
    M'_{{k} \over {2}}.
\end{array}$

 There are two cases to be considered.\\[5mm]
{\bf  Case 1.} $k=4l$. \\
[5mm] \indent In this case the induced subgraphs
$A=<u_{1'},u_{2'},...,u_{{({{k} \over {2}})}'}>$ and $B=<v_1,
v_2,...,v_{{k}\over {2}} >$  are complete graphs $K_{{k} \over
{2}}$. So they are 1-factorable. Let  $F_1, F_2,..., F_{{{k}\over{
2}}-1}$ and $F'_1, F'_2,..., F'_{{{k}\over{ 2}}-1}$ be standard
1-factorizations (see~\cite{behzad}, page 166) of $A$ and $B$,
respectively, such that the edges $u_{t'} u_{{({{k} \over {2}})
}'} \in F_t$ and $v_t v_{{k} \over {2} } \in F'_t$. Now for each
$t$ $(1 \leq t \leq {s-{{{k} \over {2}}}+1} \leq {{{k} \over
{2}}-1})$ we delete all of the edges of $F_t \backslash  \{u_{t'}
u_{{({{k} \over {2}})}' } \}$ and $F'_t$.
 Also we delete the edge $ v_{{({t+1})}'}v_{{({{k} \over {2}}) }'}$
 $\pmod {({{k} \over {2}}-1)}$.
 If there exist some edges such as
$v_i v_{i+1} \in F'_t$, then instead of these edges
we delete the edges $v_{i'} v_{i+1} \in M'_{{k} \over {2}}$.

Also for an arbitrary index $i$ of such as edges  $v_i v_{i+1}$
 we delete the edge $v_i v_{{({{k} \over {2}}) }'}$ instead of the edge
 $ v_{{({t+1})}'} v_{{({{k} \over {2}}) }'}$.
Finally we delete the edges $x_1 u_{t+1}$, $x_2 u_{t+2}$,...,
$x_{{{k} \over 2}-1} u_{t+{{k} \over 2 }-1}$ $\pmod {{{k} \over 2}}$.

We join $y_t$ to the ends of all  deleted edges.
It can be easily seen that $deg(y_t)=2(k-1)$ and the color of $y_t$ is forced
to be $t+{{k} \over {2}}$, for $t \neq {{k} \over {2}}-1$ and the color of
 $y_{{{k} \over {2}}-1}$ to be  ${{k} \over {2}}-1$.\\
 [5mm]
{\bf Case 2.} $k=4l+2$.\\
[5mm] \indent In this case the induced subgraphs
$A=<u_{1'},u_{2'},...,u_{{({{k} \over {2}})}'}, u_1>$ and
$B=<v_1, v_2,...,v_{{k}\over {2}}, v_{{({{k} \over {2}})}'}>$ are
complete graphs
 $K_{{{k}\over {2}}+1}$. So they are 1-factorable.
Let $F_1, F_2,..., F_{{k}\over{ 2}}$ and
$F'_1, F'_2,...,F'_{{k}\over{ 2}}$ be 1-factorizations of
$A$ and $B$, respectively, such that $u_1 u_{t'}  \in F_t$ and
$ v_t v_{{({{k} \over {2}}) }'} \in F'_t$.
Now for each $t$
$(1 \leq t \leq {s-{{k} \over {2}}+1} \leq {{{k} \over {2}}-1})$
we delete all of the edges of
$F_t \backslash \{u_1 u_{t'} , u_{j'} u_{{({{k} \over {2}})}'} \}$
 and $F'_t$.
Also we delete the edge  $u_j u_{j'} \in M_1$.
 If there exist some edges such as
$v_i v_{i+1} \in F'_t$ then instead of the edges $v_i v_{i+1}$
we delete the edges $v_{i'} v_{i+1} \in M'_{{k} \over {2}}$.
 Finally we delete the edges $x_1 u_{j+1}$, $x_2 u_{j+2}$,...,
$x_{{{k} \over 2}-1} u_{j+{{{k} \over 2 }}-1}$ $\pmod {{{k} \over 2}}$.
We join $y_t$ to the ends of all  deleted edges.
It can be easily seen that $deg(y_t)=2(k-1)$ and the color  of $y_t$ is forced
 to be $t+{{k} \over {2}}$.
}\end{proof}

To illustrate the construction shown in the  proof of  Theorem~\ref{keven},
we provide the following example.
 %
\begin{ex} {\rm
Let $k=8$. For $17 \leq n \leq 22$, we construct a $14$-regular
$8$-chromatic graph of order $n$  with a defining set of size $7$.
For $n=16+s$, $1 \le s \le 6$, we add $s$ new vertices to the
$14$-regular $8$-chromatic graph  $G_{2(8)}$ of order $16$ and
delete some nonessential edges as  explained in the proof of
Theorem~\ref{keven}.
\begin{table}[h]
\begin{center}{
\caption{\label{t5.2} New vertices and deleted edges.}
\vspace{.25cm}
\begin{tabular}{|c|c|c|c|c|c|c|}
\hline
New vertices  &$x_1$         & $x_2$       & $x_3$      & $y_1$     & $y_2$        & $y_3$        \\
\hline
        &$u_1u_{2'}$ & $u_1u_{3'}$ & $u_2u_{1'}$ & $u_{2'}u_{3'}$ & $u_{1'}u_{3'}$& $u_{1'}u_{2'}$\\
        &$u_2u_{3'}$ & $u_3u_{1'}$ & $u_3u_{2'}$ & $v_1v_4$       & $v_2v_4$      & $v_{3'}v_4$  \\
Deleted &$u_4u_{1'}$ & $u_4u_{2'}$ & $u_4u_{3'}$ & $v_{2'}v_{3}$ & $v_1v_3$      & $v_{1'}v_2$  \\
edges   &$v_1v_{1'}$ & $v_1v_{2'}$ & $v_1v_{3'}$ & $v_2v_{4'}$    & $v_{3'}v_{4'}$& $v_1v_{4'}$   \\
        &$v_2v_{2'}$ & $v_2v_{3'}$ & $v_3v_{1'}$ & $x_1u_2$       & $x_1u_3$      & $x_1u_4$     \\
        &$v_3v_{3'}$ & $v_4v_{1'}$ & $v_4v_{2'}$ & $x_2u_3$       & $x_2u_4$      & $x_2u_1$     \\
        &$u_3v_4$    & $u_2v_3$    & $u_1v_2$    & $x_3u_4$       & $x_3u_1$      & $x_3u_2$    \\
\hline
\end{tabular}
}\end{center}
\end{table}

Table~\ref{t5.2} gives  all  the deleted edges of $G_{2(8)}$ with
respect to addition of new vertices. In Figure~\ref{f5.3}, we
show the deleted edges and the added edges to construct a
$14$-regular $8$-chromatic graph $H_1$ of order $17$ $(s=1)$ with
a defining set of size $7$.  The dotted lines are the deleted
edges and the vertices of the defining set are shown by the
filled circles.
\hspace{-1cm}\begin{figure}[h]
 \unitlength=1mm
 \epsfysize=6cm
 \hspace{-2.5cm}\epsfbox[-35 15  500 120]{G2(8)2.eps}
\caption{\label{f5.3} $d(H_1,\chi=8)=7$.}
\end{figure}
}\end{ex}
\begin{thm}
For each $k \geq 4$,  $2k \leq n \leq 3k-2$, and $r > 2(k-1)$,
we have

$d(n , r, \Cchi=k)=k-1$.
\end{thm}
\begin{proof}{
Let $n=2k+s$, $0 \leq s \leq k-2$, and $r=2(k-1)+t$. By  Remark~1,
if there exists an $r$-regular $k$-chromatic
graph with $n$ vertices then it is obvious that $t < s$.
 We construct an $r$-regular $k$-chromatic graph $H$ with $n$ vertices
 in the following manner. \\
 Consider  graph $G_{2(k)}$,
 let $V(G_1)=\{u_1,...,u_k\}$ and  $V(G_2)=\{v_1,...,v_k\}$,
 and $c(u_i)=c(v_i)=i$, for $i=1,2,...,k$.
 We add $s$ new vertices say $x_1,...,x_s$ to $G_{2(k)}$. For each $x_i$
 $(1 \leq i \leq s)$ we join $x_i$ to each vertex of
 $V(G_1) \cup V(G_2) \backslash \{u_i, v_i\}$.
 Also, in the case of $s$ even, we add the edges of $t$
mutually disjoint 1-factors of $K_s$,
and in the case of $s$ odd, the edges of $\frac{t}{2}$  mutually
 disjoint 2-factors of $K_s$,  to $x_1,...,x_s$.
 The  graph obtained in this way, say $H'$, is a $k$-chromatic graph with $n$ vertices
 and a defining set $S=\{x_2,...,x_s,v_{s+1},...,v_k\}$ such that
 $deg(x_i)=2(k-1)+t$ $(1 \leq i \leq s)$, $deg(u_i)=deg(v_i)=2(k-1)+s-1$
 $(1 \leq i \leq s)$, and $deg(u_i)=deg(v_i)=2(k-1)+s$ $(s+1 \leq i \leq k)$.
 Now we show that by deleting some suitable nonessential edges of $H'$
  the desired $r$-regular graph $H$ can be obtained.\\
In the graph $H'$, for convenience let
 $A=\{u_1,...,u_{\lfloor{\frac{s}{2}}\rfloor}\}$,
 $C=\{u_{\lfloor{\frac{s}{2}}\rfloor+1},...,u_s\}$,
  $D=\{u_{s+1},...,u_{s+{\lfloor {\frac{k-s}{2}}\rfloor}}\}$, and
  $B=\{u_{s+{\lfloor {\frac{k-s}{2}}\rfloor} +1},...,u_k\}$.
  Also let $A'=\{v_1,...,v_{\lfloor{\frac{s}{2}}\rfloor}\}$,
 $C'=\{v_{\lfloor{\frac{s}{2}}\rfloor+1},...,v_s\}$,
  $D'=\{v_{s+1},...,v_{s+{\lfloor {\frac{k-s}{2}}\rfloor}}\}$, and
  $B'=\{v_{s+{\lfloor {\frac{k-s}{2}}\rfloor} +1},...,v_k\}$.
Let $i'=i+{\lfloor {\frac{k-s}{2}}\rfloor}$ for
$s+1 \leq i \leq s+{\lfloor {\frac{k-s}{2}}\rfloor}$.

First we delete a maximal matching of each complete bipartite subgraph
 with parts $B$ and $D$ of $G_1$ and parts $B'$ and $D'$ of $G_2$.
For $k-s$ odd, we assume $u_{k-1}$ and $v_k$ to be vertices unsaturated by the maximal matchings. Then we delete the edge $u_{k-1}v_k$.\\
Secondly,  we delete the edges of $s-t-1$ mutually disjoint maximal matchings of
 each complete bipartite subgraph  with parts
 $A \cup B$ and $C \cup D$ of $G_1$ and parts
$A' \cup B'$ and $C' \cup D'$ of $G_2$.
 For $k$ odd, we assume that the following vertices are unsaturated
by the maximal matchings: $\{u_1,...,
u_{\lfloor{\frac{s}{2}}\rfloor}, u_{(s+1)'},...,
u_{(s+1)'+s-t-2-{\lfloor{\frac{s}{2}}\rfloor}}\}$ and $\{v_2,...,
v_{\lfloor{\frac{s}{2}}\rfloor}, v_1, v_{(s+1)'+1},...\\,
v_{(s+1)'+s-t-1-{\lfloor{\frac{s}{2}}\rfloor}}\}$, in the case of
$s$ even, or
$\{u_{\lfloor{\frac{s}{2}}\rfloor+1},...,\linebreak[0]
u_{\lfloor{\frac{s}{2}}\rfloor+s-t-1}\}$ and
$\{v_{\lfloor{\frac{s}{2}}\rfloor+2},..., v_s,
v_{\lfloor{\frac{s}{2}}\rfloor+1},v_{s+2},...,
v_{2s-t-1-\lfloor{\frac{s}{2}}\rfloor}\}$, in the case of  $s$
odd. Then we delete the edges $u_1v_2, u_2v_3,...,
u_{\lfloor{\frac{s}{2}}\rfloor-1}v_{\lfloor{\frac{s}{2}}\rfloor},
u_{\lfloor{\frac{s}{2}}\rfloor}v_1, u_{(s+1)'}v_{(s+1)'+1},...,
\linebreak u_{(s+1)'+s-t-2-{\lfloor{\frac{s}{2}}\rfloor}}
v_{(s+1)'+s-t-1-{\lfloor{\frac{s}{2}}\rfloor}}$, or   the edges
$u_{\lfloor{\frac{s}{2}}\rfloor+1}v_{\lfloor{\frac{s}{2}}\rfloor+2},
..., u_sv_{\lfloor{\frac{s}{2}}\rfloor+1}\\,u_{s+1}v_{s+2},...,
\linebreak[0] u_{\lfloor{\frac{s}{2}}\rfloor+s-t-1}
v_{2s-t-1-\lfloor{\frac{s}{2}}\rfloor}$,
depending on the parity of $s$, respectively.\\
If $s-t >  \lfloor {\frac{k}{2}} \rfloor$
then in the second  step we delete $\lfloor {\frac{k}{2}} \rfloor-1$
maximal matchings.
Finally we delete the edges of $s-t-{\lfloor {\frac{k}{2}} \rfloor}$
mutually disjoint 1-factors $F_j$
$(1 \leq j \leq s-t-{\lfloor {\frac{k}{2}} \rfloor})$
of bipartite subgraph with parts $C \cup D$ and $C' \cup D'$,
where

$\begin{array}{ccll}
&F_j&=&\{u_i v_{i+j+1} \ | \ \lfloor{\frac{s}{2}}\rfloor+1 \leq i \leq
s+\lfloor{{\frac{k-s}{2}}}\rfloor-j-1\} \cup \\
&&&\{u_i v_{i-i_0+\lfloor{\frac{s}{2}}\rfloor+1} \ | \
 i_0=s+\lfloor{{\frac{k-s}{2}}}\rfloor-j \leq i \leq
 s+\lfloor{{\frac{k-s}{2}}}\rfloor\}.
\end{array}$

In fact if we consider the order
$u_{\lfloor{\frac{s}{2}}\rfloor+1},...,u_s,
u_{s+1},...,u_{s+{\lfloor {\frac{k-s}{2}}\rfloor}}$, and
$v_{\lfloor{\frac{s}{2}}\rfloor+1},...,v_s, \linebreak[0]
v_{s+1},...,v_{s+{\lfloor {\frac{k-s}{2}}\rfloor}}$,
 for the vertices in $C \cup D$ and $C' \cup D'$, respectively, then
 each  \mbox {1-factor} $F_j$ contains the edges in which  the $i$th vertex in
 $C \cup D$ is matched with $(i+j+1)$th
 vertex $\pmod {|C \cup D|}$ in $C' \cup D'$. (See Figure~4.)

 Also for decreasing the degree of vertex sets $A
 \cup B$ and $A' \cup B'$, we delete the edges of
$s-t-{\lfloor {\frac{k}{2}} \rfloor}$ mutually disjoint 1-factors
$F'_j$
 $(1 \leq j \leq s-t-{\lfloor {\frac{k}{2}} \rfloor})$
 of bipartite subgraph with parts $A \cup B$ and $A' \cup B'$  the same
 as above.
Therefore the  graph $H$ obtained in this way contains a complete
subgraph say $K_k=<A \cup B \cup C' \cup D'>$ and $H$ is an
$r$-regular graph.
\begin{center}
\begin{tabular}{rccl}
$u_{\left\lfloor{\frac{s}{2}}\right\rfloor+1}$ &\hspace{2cm} & \hspace{2cm} & $v_{\left\lfloor{\frac{s}{2}}\right\rfloor+1}$\\
$u_{\left\lfloor{\frac{s}{2}}\right\rfloor+2}$ &\hspace{2cm} & \hspace{2cm} &$v_{\left\lfloor{\frac{s}{2}}\right\rfloor+2}$ \\
&\hspace{2cm}& \hspace{2cm}  &$v_{\left\lfloor{\frac{s}{2}}\right\rfloor+3}$ \\
. &\hspace{2cm}&  \hspace{2cm} & $v_{\left\lfloor{\frac{s}{2}}\right\rfloor+4}$ \\
. &\hspace{2cm}.& \hspace{2cm} & . \\
. &\hspace{2cm}. & \hspace{2cm} & . \\
 &\hspace{2cm}. & \hspace{2cm} & . \\
$u_{s-1}$ & \hspace{2cm} & \hspace{2cm} & \\
$u_{s}$   &\hspace{2cm}  & \hspace{2cm} & $v_{s}$ \\
$u_{s+1}$ &\hspace{2cm}  & \hspace{2cm} & $v_{s+1}$ \\
. &\hspace{2cm}. & \hspace{2cm} & $v_{s+2}$\\
. &\hspace{2cm}. & \hspace{2cm} &. \\
. &\hspace{2cm}. & \hspace{2cm} &. \\
$u_{s+{\left\lfloor {\frac{k-s}{2}}\right\rfloor}-2}$ &\hspace{2cm} &  \hspace{2cm} &. \\
$u_{s+{\left\lfloor {\frac{k-s}{2}}\right\rfloor}-1}$ &\hspace{2cm} &  \hspace{2cm} & \\
$u_{s+{\left\lfloor {\frac{k-s}{2}}\right\rfloor}}$   &\hspace{2cm} &  \hspace{2cm} & $v_{s+{\left\lfloor {\frac{k-s}{2}}\right\rfloor}}$  \\
\end{tabular}
\begin{figure}[h]
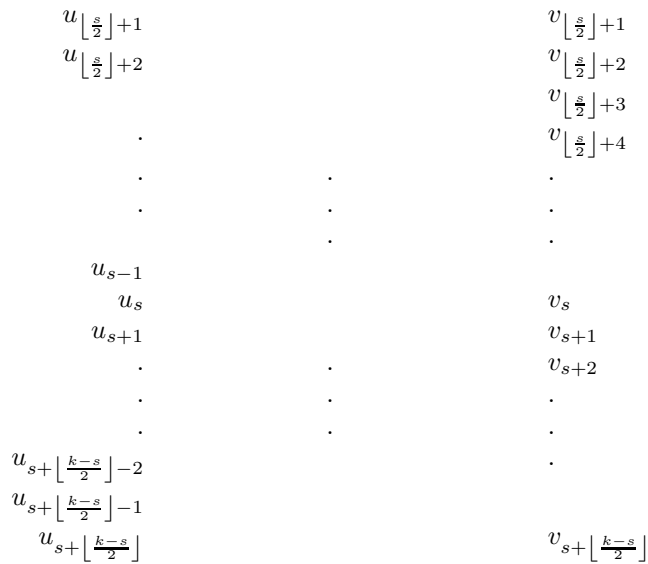

\caption{\label{f5.4} $1$-factor $F_1$.}
\end{figure}
\end{center}

\begin{figure}[h]
\unitlength=1mm
 \centering
 \epsfysize=9.5cm
 \epsfxsize=8cm
 \epsfbox[-11 -124 93 -12]{1f.eps}
\vspace{-9cm}
\end{figure}

}\end{proof}
\section*{Acknowledgments}
The authors thank Professor E.S.  Mahmoodian
for reading the manuscript and for his  helpful suggestions.


\begin{thebibliography}{1}

\bibitem{behzad}
M.~Behzad, G.~Chartrand, and L.~Lesniak.
\newblock {\em Graphs and digraphs}.
\newblock Prindle, Boston, 1979.


\bibitem{brooks}
R.L.~Brooks.
\newblock {\em On coloring the nodes of a network}.
\newblock Mathematical Proceedings of the Cambridge Philosophical Society,
{\bf 37}:194--197, 1941


\bibitem{k96}
A.D. Keedwell.
\newblock Critical sets for latin squares, graphs and block designs: a survey.
\newblock {\em Congr. Numer.}, {\bf 113}:231--245, 1996.
\newblock Festschrift for C. St. J. A. Nash-Williams.

\bibitem{kerman}
E.S. Mahmoodian.
\newblock Some problems in graph colorings.
\newblock In S.~H. Javadpour and M.~Radjabalipour, editors, {\em Proc. $26$th
  Annual Iranian Math. Conference}, pages 215--218, Kerman, March 1995. Iranian
  Math. Soc., University of Kerman.

\bibitem{mahmend}
E.S. Mahmoodian and E.~Mendelsohn.
\newblock On defining numbers of vertex coloring of regular graphs.
\newblock {\em Discrete Mathematics}, {\bf 197/198}:543--554, 1999.

\bibitem{mnz}
E.S. Mahmoodian, R.~Naserasr, and M.~Zaker.
\newblock Defining sets in vertex coloring of graphs and latin rectangles.
\newblock {\em Discrete Mathematics}, {\bf 167/168}:451--460, 1997.

\bibitem{rkgraph}
E.S.~Mahmoodian, B.~Omoomi, and N.~Soltankhah.
\newblock Smallest  defining number of $r$-regular $k$-chromatic graphs: $r\neq k$.
\newblock {\em Ars Combinatoria}, {\bf 78}: 211--223, 2006.



\bibitem{mahsol}
N.~Soltankhah and E.S.~Mahmoodian.
\newblock On defining numbers of  $k$-chromatic $k$-regular graphs.
\newblock Ars Combinatoria, {\bf 76}: 257-276, 2005.


\bibitem{west}
D.B. West.
\newblock {\em Introduction to Graph Theory}.
\newblock 2nd Eddition, Prentice Hall, Upper Saddle River, NJ, 2001.


\end{thebibliography}

\noindent E-mail addresses: \\
bomoomi@cc.iut.ac.ir \\
soltan@alzahra.ac.ir
\end{document}